\documentclass{amsart}
\usepackage{amssymb, amsthm, amsmath, amscd}
\usepackage[neveradjust]{paralist}

\newtheorem{theorem}{Theorem}[section]
\newtheorem{lemma}[theorem]{Lemma}

\newtheorem{corollary}[theorem]{Corollary}

\theoremstyle{definition}
\newtheorem{definition}[theorem]{Definition}
\newtheorem{remark}[theorem]{Remark}
\newtheorem{example}[theorem]{Example}
\newtheorem{question}[theorem]{Question}

\numberwithin{equation}{theorem}

\def\phi{\varphi}
\def\Phi{\varPhi}
\def\bar{\overline}
\def\tilde{\widetilde}
\def\dlim{\varinjlim}
\def\la{\langle}
\def\ra{\rangle}
\def\ge{\geqslant}
\def\le{\leqslant}
\def\bsf{{\boldsymbol{f}}}
\def\bsx{{\boldsymbol{x}}}

\def\pd{\operatorname{pd}}
\def\rad{\operatorname{rad}}

\def\Tor{\operatorname{Tor}}
\def\Hom{\operatorname{Hom}}
\def\Ext{\operatorname{Ext}}
\def\depth{\operatorname{depth}}
\def\to{\longrightarrow}
\def\into{\lhook\joinrel\relbar\joinrel\rightarrow}
\def\Pdot{P_\bullet}

\def\fraka{\mathfrak a}
\def\frakb{\mathfrak b}
\def\frakm{\mathfrak m}
\def\frakp{\mathfrak p}

\def\KK{\mathbb K}
\def\NN{\mathbb N}
\def\PP{\mathbb P}
\def\QQ{\mathbb Q}
\def\ZZ{\mathbb Z}

\begin{document}
\title{Local cohomology and pure morphisms}

\author{Anurag K. Singh}
\address{Department of Mathematics, University of Utah, 155 South 1400 East, Salt Lake City, UT~84112, USA} \email{singh@math.utah.edu}
\author{Uli Walther}
\address{Department of Mathematics, Purdue University, 150 N. University Street, West Lafayette, IN~47907, USA} \email{walther@math.purdue.edu}

\thanks{A.K.S.~was supported by NSF grants DMS~0300600 and DMS~0600819. U.W.~was supported by NSF grant DMS~0555319 and by NSA grant H98230-06-1-0012.}

\dedicatory{Dedicated to Professor Phil Griffith}

\begin{abstract}
We study a question raised by Eisenbud, Musta\c t\v a, and Stillman regarding the injectivity of natural maps from $\Ext$ modules to local cohomology modules. We obtain some positive answers to this question which extend earlier results of Lyubeznik. In the process, we also prove a vanishing theorem for local cohomology modules which connects theorems previously known in the case of positive characteristic and in the case of monomial ideals.
\end{abstract}
\maketitle

\section{Introduction}

Throughout this paper, the rings we consider are commutative, Noetherian, and contain an identity element. For an ideal $\fraka$ of a ring $R$, the local cohomology modules $H^i_\fraka(R)$ may be obtained as
\[
H^i_\fraka(R)=\dlim_t\Ext^i_R(R/\fraka_t,R)\qquad\text{ for }i\ge0\,,
\]
where $\{\fraka_t\}_{t\ge0}$ is a decreasing chain of ideals cofinal with the chain $\{\fraka^t\}_{t\ge0}$, and the maps in the directed system are those induced by the natural surjections
\[
R/\fraka_{t+1}\to R/\fraka_t\,.
\]
Any chain of ideals which is cofinal with the chain $\{\fraka^t\}_{t\ge0}$ yields the same direct limit. In this context, Eisenbud, Musta\c t\v a, and Stillman have raised the following questions:

\begin{question}\cite[Question~6.1]{EMS}\label{Q1}
Let $R$ be a polynomial ring over a field. For which ideals $\fraka$ of $R$ does there exist a chain of ideals $\{\fraka_t\}_{t\ge0}$ as above, such that for all $i\ge0$ and all $t\ge0$, the natural map
\[
\Ext^i_R(R/\fraka_t,R)\to H^i_\fraka(R)
\]
is injective?
\end{question}

\begin{question}\cite[Question~6.2]{EMS}\label{Q2}
Given a polynomial ring $R$ over a field, for which ideals $\fraka$ is the natural map $\Ext^i_R(R/\fraka,R)\to H^i_\fraka(R)$ an inclusion?
\end{question}

Question~\ref{Q1} is motivated by the fact that the $R$-modules $H^i_\fraka(R)$ are typically not finitely generated, whereas modules of the form $\Ext^i_R(R/\frakb,R)$ are finitely generated. Consequently, a chain of ideals as in Question~\ref{Q1} yields a filtration of $H^i_\fraka(R)$ by a natural family of finitely generated submodules.

Let $R$ be a polynomial ring over a field. For an ideal $\fraka$ generated by square-free monomials $m_1,\dots,m_r$, set $\fraka^{[t]}=(m_1^t,\dots,m_r^t)$ for integers $t\ge0$. Lyubeznik \cite[Theorem~1~(i)]{Lyubeznik-monomial} proved that the natural maps
\[
\Ext^i_R(R/\fraka^{[t]},R) \to H^i_\fraka(R)
\]
are injective for all $i\ge0$ and $t\ge0$; see also Musta\c t\v a \cite[Theorem~1.1]{Mustata}. If $R$ has positive characteristic, an ideal $\fraka$ generated by square-free monomials has the property that $R/\fraka$ is \emph{$F$-pure}; see \S\,\ref{sec:pure}. Our main result, Theorem~\ref{thm:main}, recovers Lyubeznik's result and also provides a positive answer to Question~\ref{Q1} for ideals defining $F$-pure rings, a case we single out for mention here:

\begin{theorem}\label{thm:fpure}
Let $R$ be a regular ring containing a field of characteristic $p>0$, and $\fraka$ an ideal such that $R/\fraka$ is $F$-pure. Then the natural maps
\[
\Ext^i_R(R/\fraka^{[p^t]},R)\to H^i_\fraka(R)
\]
are injective for all $i\ge 0$ and all $t\ge 0$.
\end{theorem}

\begin{remark}\label{rem:at_depth}
If $d=\depth_R(\fraka,R)$, then the natural map
\[
\Ext^d_R(R/\fraka,R)\to H^d_\fraka(R)
\]
is injective. To see this, let $E^\bullet$ be a minimal injective resolution of $R$. Then $H^\bullet_\fraka(R)$ is the cohomology of the complex $\Gamma_\fraka(E^\bullet)$ and $\Ext^\bullet_R(R/\fraka,R)$ is the cohomology of its subcomplex $\Hom_R(R/\fraka,E^\bullet)=(0:_{E^\bullet}\fraka$). Since $d$ is the least integer $i$ such that $\Gamma_\fraka(E^i)$ is nonzero, we are considering the cohomology of the rows of the diagram
\[
\begin{CD}
\cdots@>>>0@>>>\Gamma_\fraka(E^d)@>>>\Gamma_\fraka(E^{d+1})@>>>\cdots\\
@. @AAA @AAA @AAA\\
\cdots@>>>0@>>>(0:_{E^d}\fraka)@>>>(0:_{E^{d+1}})\fraka@>>>\cdots\,,
\end{CD}
\]
and the desired inclusion follows.
\end{remark}

\begin{remark}\label{rem:ci}
It is easy to see that Question~\ref{Q1} has a positive answer if $\fraka$ is a set-theoretic complete intersection: if $f_1,\dots,f_n$ is a regular sequence generating $\fraka$ up to radical, then the ideals $\fraka_t=(f_1^t,\dots,f_n^t)$ form a descending chain with $\Ext^i_R(R/\fraka_t,R)\into H^i_\fraka(R)$ for all $i\ge0$ and $t\ge0$; for $i=n$, this follows from Remark~\ref {rem:at_depth},
whereas if $i\neq n$, then $\Ext^i_R(R/\fraka_t,R)=0=H^i_\fraka(R)$.
\end{remark}

We thank David Eisenbud, Srikanth Iyengar, and Oana Veliche for useful discussions. Our work also owes a great intellectual debt to Gennady Lyubeznik's paper \cite{Lyubeznik-dual}.

\section{Pure homomorphisms and $F$-pure rings}\label{sec:pure}

\begin{definition}\label{def:pure}
A ring homomorphism $\phi\colon R\to S$ is \emph{pure} if the map
\[
\phi\otimes1\colon R\otimes_R M\to S\otimes_R M
\]
is injective for each $R$-module $M$. If $R$ contains a field of characteristic $p>0$, then $R$ is \emph{$F$-pure} if the Frobenius homomorphism $r\mapsto r^p$ is pure.
\end{definition}

Evidently, pure homomorphisms are injective. Let $R$ be a subring of $S$. If the inclusion $R\into S$ splits as a maps of $R$-modules, then it is pure. The converse is also true for module-finite extensions, see \cite[Corollary~5.3]{HR2}.

\begin{example}\label{ex:monomial}
Let $R=\KK[x_1,\dots,x_d]$ be a polynomial ring over a field $\KK$, and let $t$ be a positive integer. Then there is a $\KK$-linear endomorphism $\phi$ of $R$ with $\phi(x_i)=x_i^t$ for $1\le i\le d$. The inclusion $\phi(R)\subseteq R$ splits since $R$ is a free module over $\phi(R)$ with basis $x_1^{e_1}\cdots x_d^{e_d}$ where $0\le e_i\le t-1$. It follows that $\phi\colon R\to R$ is pure.

Let $\fraka$ be an ideal of $R$ generated by square-free monomials. Then $\phi(\fraka)\subseteq\fraka$, so $\phi$ induces an endomorphism $\bar{\phi}$ of $R/\fraka$. The image of $\bar{\phi}$ is spanned, as a $\KK$-vector space, by those monomials in $x_1^t,\dots,x_d^t$ which are not in $\fraka$. Using the map which is the identity on these monomials, and kills the rest, we obtain a splitting of $\bar{\phi}$. It follows that the endomorphism $\bar{\phi}\colon R/\fraka\to R/\fraka$ is pure.
\end{example}

\begin{remark}
The notion of $F$-pure rings was introduced by Hochster and Roberts in the course of their study of rings of invariants \cite{HR1,HR2}. Examples of $F$-pure rings include regular rings, determinantal rings, Pl\"ucker embeddings of Grassmannians, polynomial rings modulo square-free monomial ideals, normal affine semigroup rings, homogeneous coordinate rings of ordinary elliptic curves, and, more generally, homogeneous coordinate rings of ordinary Abelian varieties. Moreover, pure subrings of $F$-pure rings are $F$-pure, and if $R$ and $S$ are $F$-pure algebras over a perfect field $\KK$, then their tensor product $R\otimes_\KK S$ is also $F$-pure.
\end{remark}

\begin{remark}
Let $\phi\colon R\to S$ be a ring homomorphism. If $f\in R$, then $\phi$ localizes to give a map $R_f\to S_{\phi(f)}$. Similarly, if $\bsf$ is a sequence of elements of $R$, then $\phi$ induces a map of \v Cech complexes
\[
\check C^\bullet_\bsf(R)\to\check C^\bullet_{\phi(\bsf)}(S)\,.
\]
Setting $\fraka=(\bsf)$, we have an induced map of local cohomology groups
\[
\phi_*\colon H^i_\fraka(R)\to H^i_{\fraka S}(S)\qquad\text{ for all }i\ge0\,.
\]
Note that for $r\in R$ and $\eta\in H^i_\fraka(R)$, we have $\phi(r)\phi_*(\eta)=\phi_*(r\eta)$.

Now suppose $\phi$ is an endomorphism of $R$ with $\rad\fraka=\rad{\phi(\fraka)R}$. Then one obtains an induced \emph{action} 
\[
\phi_*\colon H^i_\fraka(R)\to H^i_{\phi(\fraka)R}(R)=H^i_\fraka(R)\,,
\]
which is an endomorphism of the underlying Abelian group. 

The archetypal example is the one where $\phi$ is the Frobenius endomorphism of a ring $R$ of prime characteristic; in this case, for all ideals $\fraka$ of $R$ and integers $i\ge0$, there is an induced action $\phi_*$ on $H^i_\fraka(R)$ known as the \emph{Frobenius action}.
\end{remark}

If $\phi\colon R\to S$ is pure, then for all ideals $\fraka$ of $R$ and all integers $i\ge0$, the induced map $\phi_*\colon H^i_\fraka(R)\to H^i_{\fraka S}(S)$ is injective, see \cite[Corollary~6.8]{HR1} or \cite[Lemma~2.1]{HR2}. In another direction, we have the following lemma, which will be a key ingredient in the proof of Theorem~\ref{thm:main}.

\begin{lemma}\label{lem:key}
Let $(R,\frakm)$ be a local ring with a pure endomorphism $\phi$ such that $\phi(\frakm)R$ is $\frakm$-primary. Then, for all $i\ge0$, the induced action
\[
\phi_*\colon H^i_\frakm(R)\to H^i_\frakm(R)
\]
is surjective up to $R$-span, i.e., $\phi_*(H^i_\frakm(R))$ generates $H^i_\frakm(R)$ as an $R$-module. 
\end{lemma}

\begin{proof}
Consider an element $\eta\in H^i_\frakm(R)$; we need to show that it belongs to the $R$-module spanned by $\phi_*(H^i_\frakm(R))$. The descending chain of $R$-modules
\[
\la\eta,\phi_*(\eta),\phi_*^2(\eta),\dots\ra\ \supseteq\
\la\phi_*(\eta),\phi_*^2(\eta),\dots\ra\ \supseteq\
\la\phi_*^2(\eta),\phi_*^3(\eta),\dots\ra
\]
stabilizes since $H^i_\frakm(R)$ is Artinian. Hence there exists $e\ge0$ such that
\begin{equation}\label{eq:dcc}
\phi_*^e(\eta)\in\la\phi_*^{e+1}(\eta),\phi_*^{e+2}(\eta),\dots\ra\,.
\end{equation}
Let $e$ be the least such integer. If $e=0$ we are done, whereas if $e\ge1$ then the $R$-module
\[
M=\frac{\la\phi_*^{e-1}(\eta),\phi_*^e(\eta),\phi_*^{e+1}(\eta),\dots\ra}{\la\phi_*^e(\eta),\phi_*^{e+1}(\eta),\dots\ra}
\]
is nonzero. But then, by the purity of $\phi$, so is its image under
\[
\phi\otimes1\colon R\otimes_RM\to R\otimes_RM\,,
\]
which contradicts (\ref{eq:dcc}).
\end{proof}

\begin{remark}\label{rem:functor}
Let $R$ be a regular ring with a flat endomorphism $\phi$. We use $R^\phi$ to denote the $R$-bimodule which has $R$ as its underlying Abelian group, the usual action of $R$ on the left, and the right $R$-action with $r'r=\phi(r)r'$ for $r\in R$ and $r'\in R^\phi$. Let $\Phi$ be the functor on the category of $R$-modules with
\[
\Phi(M)=R^\phi\otimes_RM\,,
\]
where $\Phi(M)$ is viewed as an $R$-module via the left $R$-module structure of $R^\phi$. The iteration $\Phi^t$ is the functor with
\[
\Phi^t(M)=R^\phi\otimes_R\Phi^{t-1}(M)\qquad\text{ for }t\ge1\,,
\]
where $\Phi^0$ is interpreted as the identity functor. It is easily seen that
\[
\Phi^t(M)=R^{\phi^t}\otimes_RM\,.
\]

\begin{asparaenum}
\item There is an isomorphism $\Phi(R)\cong R$ given by $r'\otimes r\mapsto r'\phi(r)$. It follows that if $M$ is a free $R$-module, then $\Phi(M)\cong M$. For a map $\alpha$ of free modules given by a matrix $(\alpha_{ij})$, the map $\Phi(\alpha)$ is given by the matrix $(\phi(\alpha_{ij}))$. Since $\phi$ is flat, $\Phi$ is an exact functor, and so it takes finite free resolutions to finite free resolutions. If $M$ and $N$ are $R$-modules, then there are natural isomorphisms
\begin{equation}\label{eq:Phicommutes}
\Phi(\Ext^i_R(M,N))\cong\Ext^i_R(\Phi(M),\Phi(N))\qquad\text{ for all }i\ge0\,.
\end{equation}
In particular, if $\fraka$ is an ideal of $R$, then \eqref{eq:Phicommutes} implies that
\[
\Phi(\Ext_R^i(R/\fraka,R))\cong\Ext_R^i(R/\phi(\fraka)R,R)\,.
\]

\item Suppose that the ideals $\{\phi^t(\fraka)R\}_{t\ge0}$ form a descending chain cofinal with the chain $\{\fraka^t\}_{t\ge0}$. Then, for each $i\ge 0$, the above isomorphism and its iterations fit into a commutative diagram
\[
\begin{CD}
\cdots@>>>\Ext^i_R(R/\phi^t(\fraka)R,R)@>>>\Ext^i_R(R/\phi^{t+1}(\fraka)R,R)@>>>\cdots\\
& & @VVV @VVV\\
\cdots@>>>\Phi^t(\Ext^i_R(R/\fraka,R))@>>>\Phi^{t+1}(\Ext^i_R(R/\fraka,R))@>>>\cdots
\end{CD}
\]
where the maps in the top row are those induced by the natural surjections $R/\phi^{t+1}(\fraka)R\to R/\phi^t(\fraka)R$, and the vertical maps are isomorphisms. Hence the bottom row has direct limit $H^i_\fraka(R)$. It follows that $H^i_\frakm(R)\cong\Phi(H^i_\frakm(R))$.

\item Assume in addition that $(R,\frakm)$ is a regular local ring of dimension $d$, and that $\phi$ is a flat local endomorphism. In this case, the dimension formula
\[
\dim R+\dim R/\phi(\frakm)R=\dim R
\]
implies that $\phi(\frakm)R$ is $\frakm$-primary. Let $E$ denote the injective hull of $R/\frakm$ as an $R$-module, and set $(-)^\vee=\Hom_R(-,E)$. Since $R$ is Gorenstein, we have
\[
E\cong H^d_\frakm(R)\cong\Phi(H^d_\frakm(R))\cong\Phi(E)\,.
\]
Hence \eqref{eq:Phicommutes} implies that $\Phi(M^\vee)\cong(\Phi(M))^\vee$ for each $R$-module $M$. Setting $M=\Ext^i_R(R/\fraka,R)$ and using local duality, we get
\[
(\Phi(\Ext^i_R(R/\fraka,R)))^\vee\cong\Phi(H^{d-i}_\frakm(R/\fraka))\,.
\]
Since $\Phi^t(-)=R^{\phi^t}\otimes_R(-)$, we immediately obtain the isomorphisms
\[
(\Phi^t(\Ext^i_R(R/\fraka,R)))^\vee\cong\Phi^t(H^{d-i}_\frakm(R/\fraka))\qquad\text{ for all }t\ge0\,.
\]
Applying $(-)^\vee$ to the diagram in (2), we get the commutative diagram
\[
\begin{CD}
\cdots@<<<H^{d-i}_\frakm(R/\phi^t(\fraka)R)@<<<H^{d-i}_\frakm(R/\phi^{t+1}(\fraka)R)@<<<\cdots\\
& & @AAA @AAA\\
\cdots@<<<\Phi^t(H^{d-i}_\frakm(R/\fraka))@<<<\Phi^{t+1}(H^{d-i}_\frakm(R/\fraka))@<<<\cdots
\end{CD}
\]
where the vertical maps are isomorphisms, and the maps in the first row are those induced by the natural surjections $R/\phi^{t+1}(\fraka)R\to R/\phi^t(\fraka)R$.
\end{asparaenum}
\end{remark}

In the archetypal example, $R$ is a regular ring containing a field of positive characteristic, and $\phi$ is the Frobenius endomorphism. In this case, $\phi$ is flat by Kunz's theorem \cite[Theorem~2.1]{Kunz}. The functor $\Phi$ is the Peskine-Szpiro functor of \cite{PS}, and the commutative diagram in Remark~\ref{rem:functor}\,(2) is precisely that obtained by Lyubeznik in \cite[Lemma~2.1]{Lyubeznik-dual}. The following is a mild generalization of \cite[Lemma~2.2]{Lyubeznik-dual}.

\begin{lemma}\label{lem:image}
Let $(R,\frakm)$ be a regular local ring with a flat local endomorphism $\phi$, and let $\fraka$ be an ideal such that $\phi(\fraka)\subseteq\fraka$. Then $\phi$ induces an endomorphism $\bar\phi$ of $R/\fraka$, and hence an action $\bar\phi_*\colon H^i_\frakm(R/\fraka)\to H^i_\frakm(R/\fraka)$. The composition
\[\CD
R^\phi\otimes_RH^i_\frakm(R/\fraka)@>\cong>>H^i_\frakm(R/\phi(\fraka)R)@>\pi>>H^i_\frakm(R/\fraka)
\endCD\]
is the map with $r'\otimes\eta\longmapsto r'\cdot\bar\phi_*(\eta)$, where $\pi$ is the map induced by the natural surjection $R/\phi(\fraka)R\to R/\fraka$.
\end{lemma}

\begin{proof}
Since $\phi(\frakm)R$ is $\frakm$-primary, if $\bsx$ is a system of parameters for $R$, then so is its image $\phi(\bsx)$. The displayed isomorphism is a consequence of the flatness of $\phi$ as we saw in Remark~\ref{rem:functor}. To analyze this isomorphism, let $\tilde{\eta}$ be a lift of $\eta\in H^i_\frakm(R/\fraka)$ to the module $\check C^i_\bsx(R/\fraka)$ of the \v Cech complex $\check C^\bullet_\bsx(R/\fraka)$. Then 
\[
\phi(\tilde{\eta})\in\check C^i_{\phi(\bsx)}(R/\phi(\fraka)R)
\]
and the image of $r'\otimes\eta$ under the isomorphism is the image of $r'\cdot\phi(\tilde{\eta})$ in $H^i_\frakm(R/\phi(\fraka)R)$. Lastly, $\pi$ maps this to $r'\cdot\bar\phi_*(\eta)\in H^i_\frakm(R/\fraka)$.
\end{proof}

We are now ready to prove the main result:

\begin{theorem}\label{thm:main}
Let $R$ be a regular ring and $\fraka$ an ideal of $R$. Suppose $R$ has a flat endomorphism $\phi$ such that $\{\phi^t(\fraka)R\}_{t\ge0}$ is a decreasing chain of ideals cofinal with $\{\fraka^t\}_{t\ge0}$, and the induced endomorphism $\bar{\phi}\colon R/\fraka\to R/\fraka$ is pure. Then, for all $i\ge 0$ and $t\ge 0$, the natural map
\[
\Ext^i_R(R/\phi^t(\fraka)R,R)\to\Ext^i_R(R/\phi^{t+1}(\fraka)R,R)
\]
is injective.
\end{theorem}

\begin{proof}
It suffices to verify the injectivity after localizing at maximal ideals, so we assume that $(R,\frakm)$ is a regular local ring. Let $d=\dim R$, and let $E$ be the injective hull of $R/\frakm$ as an $R$-module. Using $(-)^\vee=\Hom_R(-,E)$, local duality gives an isomorphism
\[
\Ext^i_R(R/\phi^t(\fraka)R,R)^\vee\cong H^{d-i}_\frakm(R/\phi^t(\fraka)R)\,,
\]
and it suffices to show that the map
\begin{equation}\label{eq:lc-onto}
H^{d-i}_\frakm(R/\phi^{t+1}(\fraka)R)\to H^{d-i}_\frakm(R/\phi^t(\fraka)R)
\end{equation}
induced by the natural surjection
\[
R/\phi^{t+1}(\fraka)R\to R/\phi^t(\fraka)R
\]
is surjective for each $t\ge0$. In view of the isomorphisms
\[
R^\phi\otimes_R H^{d-i}_\frakm(R/\phi^t(\fraka)R)\cong H^{d-i}_\frakm(R/\phi^{t+1}(\fraka)R)
\]
and the right exactness of tensor, it suffices to verify the surjectivity of \eqref{eq:lc-onto} in the case $t=0$. By Lemma~\ref{lem:image}, this reduces to checking that the $\bar\phi$-action
\[
\bar\phi_*\colon H^{d-i}_\frakm(R/\fraka)\to H^{d-i}_\frakm(R/\fraka)\,
\]
is surjective up to taking the $R$-span. This follows from Lemma~\ref{lem:key}.
\end{proof}

Theorem~\ref{thm:fpure} follows immediately from Theorem~\ref{thm:main} by taking $\phi$ to be the Frobenius endomorphism. To recover the result for square-free monomial ideals \cite[Theorem~1~(i)]{Lyubeznik-monomial}, take $\phi$ as in Example~\ref{ex:monomial}.

\section{Examples}

We first construct an example of a module $M$ over a regular local ring $(R,\frakm)$ such that $H^i_\frakm(M)=0$ but $\Ext^i_R(R/\fraka,M)$ is nonzero for every $\frakm$-primary ideal $\fraka$ of $R$. It then follows that $H^i_\frakm(M)$ cannot be realized as a union of appropriate $\Ext$-modules. We use the following lemma:

\begin{lemma}\label{lemma:dual}
Let $(R,\frakm)$ be a regular local ring of dimension $d$, and let $\fraka$ be an $\frakm$-primary ideal. Then, for each $R$-module $M$, there is an isomorphism
\[
\Ext^i_R(R/\fraka,M)\cong\Tor^R_{d-i}(\Ext^d_R(R/\fraka,R),M)\qquad\text{ for all }0\le i\le d\,.
\]
\end{lemma}

\begin{proof}
Let $\Pdot$ be a minimal free resolution of $R/\fraka$. The complex $\Hom_R(\Pdot,R)$ has homology $\Ext^\bullet_R(R/\fraka,R)$. Since $\fraka$ is an $\frakm$-primary ideal of a regular ring $R$, we have $\depth_\fraka R=d$, and so $\Ext^j_R(R/\fraka,R)$ is nonzero only for $j=d$. It follows that, with a change of index, $\Hom_R(\Pdot,R)$ is an acyclic complex of free modules resolving the module $\Ext^d_R(R/\fraka,R)$. Hence
\begin{align*}
\Ext^i_R(R/\fraka,M)=H^i(\Hom(\Pdot,M))&\cong H^i(\Hom(\Pdot,R)\otimes_RM)\\
&\cong\Tor^R_{d-i}(\Ext^d_R(R/\fraka,R),M)\,.\qedhere
\end{align*}
\end{proof}

\begin{example}
Let $(R,\frakm)$ be a regular local ring of dimension $d>0$, and $x$ a nonzero element of $\frakm$. Then $R/(x)$ has dimension $d-1$, so $H^d_\frakm(R/(x))=0$. However, if $\fraka$ is an $\frakm$-primary ideal, then Lemma~\ref{lemma:dual} implies that
\[
\Ext^d_R(R/\fraka,R/(x))\cong\Ext^d_R(R/\fraka,R)\otimes_R R/(x)\,,
\]
which is nonzero. In particular, if $\{\fraka_t\}_{t\ge0}$ is a decreasing family of ideals cofinal with $\{\frakm^t\}_{t\ge0}$, then the modules $\Ext^d_R(R/\fraka_t,R/(x))$ are nonzero for each $t$, and so the maps $\Ext^d_R(R/\fraka_t,R/(x))\to H^d_\frakm(R/(x))$ are not injective.
\end{example}

Example~\ref{ex:de} below is due to Eisenbud: given positive integers $a\le b-2$, there exists a polynomial ring $R$ and a finitely generated graded $R$-module $M$, such that the natural map $\Ext^i_R(R/\fraka,M)\to H^i_\frakm(M)$ is not injective for all $a<i<b$ and all $\frakm$-primary ideals $\fraka$. This is based on a construction of Evans and Griffith \cite[Theorem~A]{EG}.

\begin{theorem}[Evans-Griffith] \label{thm:eg}
Let $\KK$ be an infinite field and take a sequence of positive integers, $n_0<n_1<\dots<n_s$. Then there exists a polynomial ring $R$ over $\KK$, with a homogeneous prime ideal $\frakp$, such that the local cohomology module $H^i_\frakm(R/\frakp)$ is nonzero if and only if $i\in\{n_0,n_1,\dots,n_s\}$. Moreover, if $n_0\ge 2$, then $R/\frakp$ may be chosen to be a normal domain.\qed
\end{theorem}

\begin{example}[Eisenbud]\label{ex:de}
Let $a\le b-2$ be positive integers. By Theorem~\ref{thm:eg}, there exists a polynomial ring $R$ with a homogeneous prime $\frakp$, such that $\depth R/\frakp=a$, $\dim R/\frakp=b$ and $H^j_\frakm(R/\frakp)=0$ for all $a<j<b$. Let $\fraka$ be an $\frakm$-primary ideal. Then $\Ext^a_R(R/\fraka,R/\frakp)$ is nonzero so, by Lemma~\ref{lemma:dual},
\[
\Tor^R_{d-a}(\Ext^d_R(R/\fraka,R),R/\frakp)\neq 0\qquad\text{ where }d=\dim R\,.
\]
By the rigidity of $\Tor$ over regular local rings, \cite{rigidity}, it follows that
\[
\Tor^R_j(\Ext^d_R(R/\fraka,R),R/\frakp)\neq 0\qquad\text{ for all }0\le j\le d-a\,.
\]
By another application of Lemma~\ref{lemma:dual}, the module $\Ext^i_R(R/\fraka,R/\frakp)$ is nonzero if $a\le i\le d$. Now if $\{\fraka_t\}_{t\ge0}$ is any decreasing family of ideals cofinal with $\{\frakm^t\}_{t\ge0}$, it follows that the maps
\[
\Ext^i_R(R/\fraka_t,R/\frakp)\to H^i_\frakm(R/\frakp)
\]
are not injective for each $a<i<b$ and each $t\ge0$.
\end{example}

\begin{example}\label{ex:0134}
Let $\KK$ be a field and consider the $\KK$-linear ring homomorphism
\[
\alpha\colon R=\KK[w,x,y,z]\to\KK[s^4,s^3t,st^3,t^4]
\]
where $\alpha$ sends $w,x,y,z$ to the elements $s^4,s^3t,st^3,t^4$ respectively. Let $\fraka$ be the kernel of $\alpha$. Using vanishing theorems such as \cite[Theorem~2.9]{HL}, it may be verified that $H^i_\fraka(R)=0$ for $i\ge 3$. 

If $\KK$ has characteristic $p>0$, Hartshorne \cite{Hartshorne-0134} showed that $\fraka$ is a set-theoretic complete intersection, i.e., that there exist elements $f,g$ in $R$ such that $\fraka=\rad(f,g)$. In this case, the ideals $\fraka_t=(f^t,g^t)$ form a descending chain cofinal with $\{\fraka^t\}$ for which the maps $\Ext^i_R(R/\fraka_t,R)\to H^i_\fraka(R)$ are injective for all $i\ge 0$ and $t\ge 0$; see Remark~\ref{rem:ci}.

Next, suppose that $\KK$ has characteristic $0$. If $\frakb$ is an ideal with $\rad\frakb=\fraka$ such that $\Ext^i_R(R/\frakb,R)\to H^i_\fraka(R)$ is injective for all $i\ge 0$, then
\[
\Ext^3_R(R/\frakb,R)=0=\Ext^4_R(R/\frakb,R)
\]
and so $R/\frakb$ is Cohen-Macaulay. This leads to the following question:
\end{example}

\begin{question}
Let $\KK$ be a field of characteristic $0$ and, as in Example~\ref{ex:0134}, let $\fraka\subset R=\KK[w,x,y,z]$ be an ideal with $R/\fraka\cong\KK[s^4,s^3t,st^3,t^4]$. Is the ideal $\fraka$ set-theoretically Cohen-Macaulay, i.e., does there exist an ideal $\frakb\subset R$ with $\rad\frakb=\fraka$, such that the ring $R/\frakb$ is Cohen-Macaulay?
\end{question}

While the requirement of $F$-purity in Theorem~\ref{thm:fpure} is certainly a strong hypothesis, it appears to be a crucial ingredient. In the following example, we have regular rings $R_p=R/pR$ of prime characteristic $p$ and ideals $\fraka_p=\fraka R_p$ such that the maps 
\[
\Ext^4_{R_p}(R_p/\fraka_p^{[p^t]},R_p)\to H^4_{\fraka_p}(R_p)
\]
are injective if and only if $R_p/\fraka_p$ is $F$-pure; the set of primes for which this is the case is infinite, as is its complement.

\begin{example}
Let $E\subset\PP^2_\QQ$ be an elliptic curve, and consider the Segre embedding of $E\times\PP^1_\QQ$ in $\PP^5_\QQ$. Clearing denominators in a set of generators for the defining ideal of the homogeneous coordinate ring, we obtain an ideal $\fraka$ of $R=\ZZ[u,v,w,x,y,z]$ such that $R/\fraka\otimes_\ZZ\QQ$ is the coordinate ring of $E\times\PP^1_\QQ$. For prime integers $p$, let $R_p=R/pR$ and $\fraka_p=\fraka R_p$. For all but finitely many primes $p$, the reduction mod $p$ of $E$ is a smooth elliptic curve $E_p$ and $R_p/\fraka_p$ is a homogeneous coordinate ring for $E_p\times\PP^1_{\ZZ/p}$. We restrict our attention to such primes. Since $\depth R_p/\fraka_p=2$, the Auslander-Buchsbaum formula implies that $\pd_{R_p}R_p/\fraka_p=4$. Using the flatness of Frobenius, we see that
\[
\pd_{R_p}R_p/\fraka_p^{[p^t]}=4\,,
\]
and hence that
\[
\Ext^4_{R_p}(R_p/\fraka_p^{[p^t]},R_p)\neq 0\qquad\text{ for all }t\ge0\,.
\]
On the other hand, $H^4_{\fraka_p}(R_p)$ is zero if $E_p$ is supersingular and nonzero if $E_p$ is ordinary, see \cite[Example~3,~page~75]{HS} or \cite[page~219]{Lyubeznik-dual}. By well-know results on elliptic curves, there are infinitely primes $p$ for which $E_p$ is supersingular, and infinitely many for which it is ordinary. Consider the natural map
\begin{equation}\label{eq:vary}
\Ext^i_{R_p}(R_p/\fraka_p^{[p^t]},R_p)\to H^i_{\fraka_p}(R_p)\,.
\end{equation}

\noindent\emph{Ordinary primes.} If $E_p$ is ordinary then its coordinate ring is $F$-pure, and it follows that $R_p/\fraka_p$ is $F$-pure as well. In this case, Theorem~\ref{thm:fpure} implies that the map \eqref{eq:vary} is injective for all $i\ge 0$ and $t\ge 0$.

\noindent\emph{Supersingular primes.} If $p$ is a prime such that $E_p$ is supersingular, then $H^4_{\fraka_p}(R_p)=0$ so the map \eqref{eq:vary} is not injective for $i=4$. We do not know whether there exists an $\fraka_p$-primary ideal $\frakb$ for which the maps
\[
\Ext^i_{R_p}(R_p/\frakb,R_p)\to H^i_{\fraka_p}(R_p)
\]
are injective for all $i\ge 0$. Since $H^i_{\fraka_p}(R_p)=0$ for $i\ge 4$ in the supersingular case, the existence of such an ideal would imply that $\fraka_p$ is set-theoretically Cohen-Macaulay; see also \cite[\S\,3]{SW}.
\end{example}

\section{A vanishing criterion}

The observations from \S\,\ref{sec:pure} yield the following vanishing theorem, which links Lyubeznik's positive characteristic result \cite[Theorem~1.1]{Lyubeznik-dual} to a theorem for monomial ideals recorded below as Corollary~\ref{cor:monomial}.

\begin{theorem}\label{thm:dual}
Let $(R,\frakm)$ be a regular local ring, $\fraka$ an ideal, and $\phi$ a flat local endomorphism such that $\{\phi^t(\fraka)R\}_{t\ge0}$ is a decreasing chain of ideals cofinal with the chain $\{\fraka^t\}_{t\ge0}$. Then $H^i_\fraka(R)=0$ if and only if some iteration of the induced action
\[
\bar\phi_*\colon H^{\dim R-i}_\frakm(R/\fraka)\to H^{\dim R-i}_\frakm(R/\fraka)
\]
is zero.
\end{theorem}

\begin{proof}
Let $d=\dim R$. The direct limit
\[
H^i_\fraka(R)=\dlim_t\Ext^i_R(R/\phi^t(\fraka)R,R)
\]
vanishes if and only if for each $t\in\NN$, there exists $k\in\NN$ such that the map 
\begin{equation}\label{eq:extmap}
\Ext^i_R(R/\phi^t(\fraka)R,R)\to\Ext^i_R(R/\phi^{t+k}(\fraka)R,R)
\end{equation}
induced by the surjection $R/\phi^{t+k}(\fraka)R\to R/\phi^t(\fraka)R$ is zero. By local duality, the map \eqref{eq:extmap} is zero if and only if
\[
H^{d-i}_\frakm(R/\phi^{t+k}(\fraka)R)\to H^{d-i}_\frakm(R/\phi^t(\fraka)R)
\]
is the zero map. By Remark~\ref{rem:functor}\,(3) and the flatness of $R^\phi\otimes_R-$, this is equivalent to the map
\[
H^{d-i}_\frakm(R/\phi^k(\fraka)R)\to H^{d-i}_\frakm(R/\fraka)
\]
being zero. By Lemma~\ref{lem:image}, this last condition is equivalent to the $k$-th iterate of the action $\bar\phi_*\colon H^{d-i}_\frakm(R/\fraka)\to H^{d-i}_\frakm(R/\fraka)$ being zero.
\end{proof}

Using Theorem~\ref{thm:dual} we next recover a vanishing theorem for monomial ideals, \cite[Theorem~1~(iii)]{Lyubeznik-monomial}. In \cite[Corollary~6.7]{Ezra} Miller proves a stronger statement connecting $H^i_\fraka(S)$ and $H^{\dim S-i}_\frakm(S/\fraka)$ via Alexander duality.

\begin{corollary}\label{cor:monomial}
Let $S$ be a polynomial ring over a field, and let $\fraka$ be an ideal generated by square-free monomials. Then $H^i_\fraka(S)=0$ if and only if $H^{\dim S-i}_\frakm(S/\fraka)=0$.
\end{corollary}

\begin{proof}
Let $S=\KK[x_1,\dots,x_d]$, and let $\phi$ be the $\KK$-linear endomorphism with $\phi(x_i)=x_i^2$ for $1\le i\le d$. Then $\phi$ is flat, and induces a pure endomorphism of $S/\fraka$, see Example~\ref{ex:monomial}.

Each of the modules in question is graded, so the issue of vanishing is unchanged under localization at the homogeneous maximal ideal of $S$. We can therefore work over the regular local ring $(R,\frakm)$, where we need to show that $H^i_\fraka(R)=0$ if and only if $H^{d-i}_\frakm(R/\fraka R)=0$. The endomorphism $\phi$ localizes to give a flat endomorphism of $R$. Moreover, since purity localizes, $\phi$ induces a pure endomorphism $\bar\phi$ of $R/\fraka R$. By Theorem~\ref{thm:dual}, $H^i_\fraka(R)=0$ if and only if some iterate of the action 
\[
\bar\phi_*\colon H^{d-i}_\frakm(R/\fraka R)\to H^{d-i}_\frakm(R/\fraka R)
\]
is zero. But $\bar\phi_*$ is injective since $\bar\phi$ is pure, so an iterate of $\bar\phi_*$ is zero precisely if $H^{d-i}_\frakm(R/\fraka R)=0$.
\end{proof}

\end{document}